\documentclass[a4paper]{article}

\usepackage{delarray,verbatim,enumerate}
\usepackage{amsmath,amsthm,amstext,amsbsy,amssymb,amsfonts,amscd}
\usepackage[latin1]{inputenc}
\usepackage[francais]{babel}
\usepackage[T1]{fontenc}

\usepackage{color}
\definecolor{vert}{rgb}{0.1,0.5,0.2}
\usepackage[colorlinks,linkcolor=blue,urlcolor=blue,citecolor=vert,linktocpage]{hyperref}

\usepackage[all]{xy}%

\newtheorem{Th}{Théorème}%[section]
\newtheorem{Lem}[Th]{Lemme}
\newtheorem{Prop}[Th]{Proposition}

\newtheorem{Sco}[Th]{Scolie}

\def\Preuve{\smallskip\noindent {\it Preuve.~}}

\def\Remarque{\smallskip\noindent {\it Remarque.~}}

\def\ie{{\it i.e. }}	\def\cf{{\it cf. }}	\def\eg{{\it e.g. }}

\def\Z{\mathbb {Z}}     			\def\Q{\mathbb Q}
\def\Zl{\mathbb{Z}_\ell}

\def\J{\mathcal  J}  				\def\R{\mathcal  R}
\def\E{\mathcal  E} 				\def\U{\mathcal  U}
				
\def\Dl{\mathcal  D\ell} 			  	\def\Cl{\mathcal  C\ell}

\def\p{\mathfrak p}					

\def\wi{\widetilde}		
	\def\deg{\operatorname{deg}}		\def\Tr{\operatorname{Tr}}
\def\Gal{\operatorname{Gal}}	
\def\Ker{\operatorname{Ker}}	
 \def\Ver{\operatorname{Ver}}

\date{}

\title{\LARGE{Classes logarithmiques et capitulation}}

\author{{\small par}\\
\\
Jean-François {\sc Jaulent}}

\begin{document}
\maketitle
\medskip

{\footnotesize
\noindent{\bf Résumé.} Nous étudions un analogue logarithmique du Théorème d'Artin-Furtwängler sur la capitulation en transposant dans le cadre des classes logarithmiques les arguments mis en \oe uvre dans la preuve algébrique classique du Théorème de l'idéal principal.}\smallskip

{\footnotesize
\noindent{\bf Abstract.} We study a logarithmic version of the classical result of Artin-Furwängler on principalization of ideal classes in the Hilbert class-field by applying the group theoretic description of the transfert map to logarithmic class-groups of degree 0.}
\smallskip\bigskip\bigskip

%\section{Introduction}
\noindent{\bf 1. Introduction}\medskip

Le célèbre Théorème d'Artin-Furwängler (\cf \cite{Art,Fur,A-T}) affirme que le groupe des classes d'idéaux $Cl_K$ d'un corps de nombres $K$ capitule dans son corps des classes de Hilbert $L=H_K$; en d'autres termes que les idéaux de $K$ se principalisent dans son extension abélienne non ramifiée maximale $L$.\smallskip

Et, comme l'extension des classes d'idéaux est injective pour les $\ell$-parties dans une extension de degré étranger à $\ell$, cela revient à dire que, pour tout nombre premier $\ell$, le $\ell$-sous-groupe de Sylow du groupe des classes d'idéaux d'un corps de nombres capitule dans le $\ell$-corps des classes de Hilbert de ce corps \ie dans sa $\ell$-extension abélienne non ramifiée maximale.\smallskip

La preuve du Théorème d'Artin-Furwängler repose sur deux éléments, le premier de nature arithmétique, le second purement algébrique: d'une part, la Théorie du corps de classes interprète le groupe de classes d'idéaux $Cl_K$ comme groupe de Galois $G_K$ de l'extension abélienne $H_K/K$, l'homomorphisme d'extension $j_{L/K}:\;Cl_K\rightarrow Cl_L$ correspondant dans cette description au morphisme de transfert $\Ver_{L/K}:\;G_K\rightarrow G_L$; d'autre part, des considérations de théorie des groupes montrent alors que le transfert est nul dans la situation étudiée:

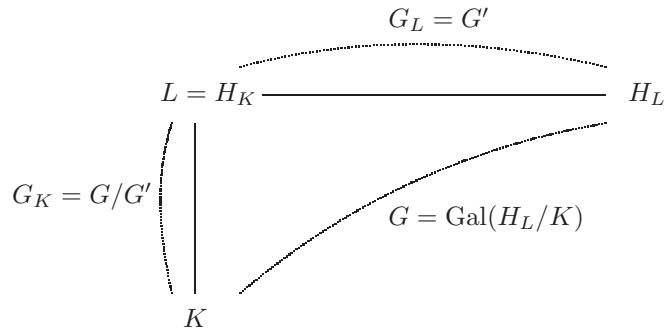
\begin{center}
\unitlength=1.5cm
\begin{picture}(6.6,3)

\put(1.7,0){$K$}
\put(1.8,0.3){\line(0,1){1.5}}
\put(1.5,2){$L=H_K$}

\bezier{60}(1.6,0.3)(1.4,1.2)(1.6,1.8)
\put(0.2,1.1){$G_K=G/G'$}

\put(2.4,2.05){\line(1,0){3.0}}
\put(5.6,2){$H_L$}

\bezier{100}(2.2,2.3)(3.7,2.7)(5.4,2.3)
\put(3.5,2.65){$G_L=G'$}

\bezier{140}(2.2,0.3)(3.6,1.5)(5.4,1.8)
\put(3.5,0.9){$G=\Gal(H_L/K)$}

\end{picture}
\end{center}
\medskip

Comme noté dans \cite{Her}, ce résultat s'étend à diverses situations arithmétiques, notamment aux groupes de classes de rayons, qui conduisent à des situations formellement comparables. Il était donc tentant d'essayer de le transposer à des groupes de classes présentant des analogies remarquables avec les objets précédents.\medskip

Il est ainsi défini dans \cite{J28}, pour chaque nombre premier $\ell$ et tout corps de nombres $K$, un groupe de classes logarithmiques noté $\,\Cl_K$, obtenu en remplaçant les valuations ordinaires $\nu_\p$ attachées à chaque place finie $\p$ de $K$ par leurs analogues formels $\wi\nu_\p$ définis à partir du logarithme $\ell$-adique de la valeur absolue $\ell$-adique attachée à la place $\p$. La Théorie $\ell$-adique du corps de classes (\cf \cite{J31, Gra}) interprète alors le groupe de classes $\,\Cl_K$ comme groupe de Galois de la pro-$\ell$-extension abélienne localement cyclotomique maximale $K^{lc}$ du corps $K$. \smallskip

Mais une première difficulté apparaît alors: la pro-$\ell$-extension abélienne $K^{lc}$ contient évidemment la $\Zl$-extension cyclotomique $K^c$ de $K$, de sorte que le groupe $\Cl_K\simeq\Gal(K^{lc}/K)$ n'est jamais fini. Pour obtenir un groupe fini, il est ainsi nécessaire de se restreindre au sous-groupe des classes logarithmiques de degré nul $\,\wi\Cl_K \simeq \Gal(K^{lc}/K^c)$. Et la {\em conjecture de Gross-Kuz'min}, souvent appelée {\em conjecture de Gross généralisée} (\cf \cite{J10,J28,J31,Kuz})-- affirme que cette restriction nécessaire est aussi suffisante, \ie que le groupe obtenu $\,\wi\Cl_K$ est effectivement un $\ell$-groupe fini.\smallskip

Lorsqu'elle est satisfaite, l'extension localement cyclotomique $K^{lc}$ provient alors, par composition avec $K^c$, d'une extension abélienne $L$ de $K$ (que l'on peut supposer linéairement disjointe de $K^c$) et il est naturel d'introduire la pro-$\ell$-extension abélienne localement cyclotomique maximale $L^{lc}$ du corps $L$ et de considérer le morphisme d'extension $\wi j_{L/K}: \;\wi\Cl_K \rightarrow \Cl_L$. Dans le schéma obtenu cependant arrive une deuxième difficulté:\medskip

\begin{center}
\unitlength=1.5cm
\begin{picture}(6.6,2.8)

\put(0.7,0){$K$}
\put(0.8,0.3){\line(0,1){1.5}}
\put(0.6,2){$K^c$}
\put(3.2,0){$L$}
\put(1.2,0.05){\line(1,0){1.6}}

\bezier{60}(0.6,0.3)(0.4,1.2)(0.6,1.8)
\put(0.3,1.0){$\Gamma$}

\put(1.2,2.05){\line(1,0){1.5}}
\put(2.9,2){$K^{lc}=L^c$}
\put(3.9,2.05){\line(1,0){1.5}}
\put(5.6,2){$L^{lc}$}
\put(3.3,0.3){\line(0,1){1.5}}

\bezier{100}(1.2,2.3)(3,2.7)(5.4,2.3)
\put(3.2,2.65){$G$}

\bezier{50}(1.2,1.8)(2.0,1.6)(2.9,1.8)
\put(1.9,1.4){$\wi\Cl_K$}

\bezier{50}(3.7,1.8)(4.5,1.6)(5.4,1.8)
\put(4.4,1.4){$\wi\Cl_L$}

\end{picture}
\end{center}
\medskip\medskip

\noindent le sous-corps $K^{lc}=L^c$ de $L^{lc}$ n'est plus alors le sous-corps maximal de $L^{lc}$ qui est abélien sur $K^c$, mais celui qui est abélien sur $K$; de sorte qu'on ne peut plus écrire comme plus haut $\wi\Cl_L \simeq G'$ et $\wi\Cl_K \simeq G/G'$ avec $G=\Gal(L^{lc}/K^c)$, mais qu'il nous faut prendre en compte aussi l'action du groupe procyclique $\Gamma=\Gal(K^c/K)$.\medskip

Le but de cette note est d'étudier la transposition dans ce contexte logarithmique des résultats classiques d'Artin-Furtwängler et de Tannaka-Terada sur la capitulation dans le corps des classes de Hilbert ou dans le corps des genres relatif à une extension cyclique non ramifiée (devenant ici procyclique).%:\medskip
\bigskip

%%%%%%%%%%%%%%%%%%%%%%%%%%%%%%%%%%%%%%%%%%%%%%%%%%%%%%%
\noindent{\bf 2. Bref rappel sur les classes logarithmiques}\medskip
%%%%%%%%%%%%%%%%%%%%%%%%%%%%%%%%%%%%%%%%%%%%%%%%%%%%%%%

Classiquement, le groupe des classes d'idéaux d'un corps de nombres $K$ est défini comme conoyau $Cl_K$ du morphisme naturel partant du groupe multiplicatif $K^\times$ à valeurs dans le groupe des idéaux $Id_K$ donné par la famille des valuations $\nu=(\nu_\p)_{\p\in Pl_K}$ attachées aux places finies de $K$.\smallskip

\centerline{$1 \rightarrow E_K \rightarrow  K^\times \overset{\nu}{\longrightarrow}  Id_K \rightarrow Cl_K \rightarrow 1$.}\smallskip

\noindent Par produit tensoriel avec $\Zl$, son $\ell$-sous-groupe de Sylow apparaît comme conoyau du morphisme $\nu$ étendu au tensorisé $\R_K=\Zl\otimes_\Z K^\times$ et à valeurs dans le $\Zl$-module libre construit sur ces mêmes places: $\Dl=\oplus_{\p\in Pl_K}\Zl\p$.
\smallskip

Le groupe des classes logarithmiques est le groupe analogue $\,\Cl_K$ obtenu en remplaçant les valuations classiques $\nu_\p$ par leurs homologues $\ell$-adiques $\wi\nu_\p$ définis à partir des logarithmes des valeurs absolues $\ell$-adiques (\cf \cite{J28}):\smallskip

\centerline{$1 \rightarrow \E_K \rightarrow  \R_K\overset{\wi\nu}{\longrightarrow}  \Dl_K \rightarrow \Cl_K \rightarrow 1$.}\smallskip

\noindent Contrairement au groupe de classes d'idéaux, c'est donc un objet $\ell$-adique.\par

Pour chaque place finie $\p$ de $K$, soit $\R_\p=\varprojlim K_\p^\times/K_\p^{\times\ell^n}$ le compactifié $\ell$-adique du groupe $K_\p^\times$ et $\J_K= \prod_\p^{r\!e\!s} \R_\p$ le $\ell$-adifié du groupe des idèles de $K$.\smallskip

 Du point de vue local, le noyau $\U_\p$ de $\nu_\p$ dans $\R_\p$ (autrement dit le sous-groupe des unités de $\R_\p$) est le groupe de normes associé à la $\Zl$-extension non ramifiée de $K_\p$; tandis que le noyau $\wi\U_\p$ de $\wi\nu_\p$ (\ie le sous-groupe des unités logarithmiques) correspond, lui, à sa $\Zl$-extension cyclotomique. Par la Théorie $\ell$-adique du corps de classes (\cf \cite{Gra,J28,J31}), le $\ell$-goupe des classes d'idéaux s'interprète comme groupe de Galois de la $\ell$-extension abélienne non ramifiée maximale $K^{nr}$ de $K$; et le $\ell$-groupe des classes logarithmiques $\,\Cl_K\simeq \J_K/\prod_\p\wi\U_\p\R_K$ comme groupe de Galois de sa pro-$\ell$-extension abélienne localement cyclotomique maximale $K^{lc}$. Le corps $K^{lc}$ est ainsi la plus grande pro-$\ell$-extension abélienne de $K$ qui est complètement décomposée %(en toutes ses places) 
 au-dessus de la $\Zl$-extension cyclotomique $K^c$. En particulier $K^{lc}$ contient $K^c$ et $\,\Cl_K$ n'est jamais fini.\smallskip

La surjection canonique du $\ell$-adifié $\J_K$ du groupe des idèles de $K$ dans le groupe procycliqe $\Gal(K^c/K)\simeq\Zl$ fournit cependant un morphisme {\em degré}:\smallskip

\centerline{$\deg:\J_K\rightarrow \Zl$;}

\noindent dont le noyau $\wi\J_K$ est, par construction, le sous-groupe normique de $\J_K$ attaché à $K^c$. Le noyau

\centerline{$\,\wi\Cl_K \simeq \wi\J_K/\prod_\p\wi\U_\p\R_K \simeq \Gal(K^{lc}/K^c)$}\medskip

\noindent de l'application induite sur $\,\Cl_K$ est ainsi le {\em groupe des classes logarithmiques de degré nul}, qui est l'objet de cette note.\smallskip

La {\em Conjecture de Gross-Kuz'min} (pour le corps $K$ et le premier $\ell$) postule précisément la finitude de ce groupe. Comme expliqué dans \cite{J10},  c'est une conséquence d'une conjecture plus générale d'indépendance $\ell$-adique de nombres algébriques, qui résulte elle-même de la conjecture de Schanuel $\ell$-adique.\smallskip

Comme établi par Greenberg, les résultats de transcendance de Baker-Brumer (\cf \eg \cite{Gra,Gre,J10}) assurent que la conjecture de Gross-Kuz'min vaut en particulier dès que le corps considéré $K$ est abélien sur $\Q$. La même conclusion vaut encore sous des hypothèses plus larges, lorsque $K$ est totalement réel, en vertu d'un raffinement du théorème de Baker-Brumer dû à M. Waldschmidt et précisé par M. Laurent puis D. Roy (\cf \cite{J38}). Le cas général reste actuellement ouvert.
\bigskip

%%%%%%%%%%%%%%%%%%%%%%%%%%%%%%%%%%%%%%%%%%%%%%%%%%%%%%%%%%%%
\noindent{\bf 3. Le schéma galoisien de la capitulation}\medskip
%%%%%%%%%%%%%%%%%%%%%%%%%%%%%%%%%%%%%%%%%%%%%%%%%%%%%%%%%%%%

Partons donc d'un corps de nombre arbitraire, vérifiant la conjecture de Gross-Kuz'min pour un premier donné $\ell$. Notons $K^c$ sa $\Zl$-extension cyclotomique de $K$ et $K^{lc}$ la pro-$\ell$-extension abélienne localement cyclotomique de $K$. Comme expliqué plus haut, le groupe $\Gal (K^{lc}/K)$ s'identifie alors au pro-$\ell$-groupe des classes logarithmiques $\Cl_K$ de $K$ et le sous-groupe $\Gal(K^{lc}/K^c)$ au $\ell$-groupe $\wi\Cl_K$ des classes logarithmiques de degré nul.

Soit maintenant $L$ une $\ell$-extension abélienne de $K$ telle qu'on ait $LK^c=K^{lc}$ et $\wi\Cl_L$ son (pro)-$\ell$-groupe des classes logarithmiques de degré nul. Il s'agit de voir que $\wi\Cl_K$ capitule dans $\wi\Cl_L$. Quitte à remplacer $L$ par un sous-corps convenable $L'$, nous pouvons supposer $L/K$ linéairement disjointe de $K^c/K$ sans restreindre aucunement la généralité de notre preuve, puisque la capitulation dans $\wi\Cl_{L'}$ entraînera la capitulation dans $\wi\Cl_L$.\smallskip

Cela fait, introduisons la pro-$\ell$-extension abélienne localement cyclotomique maximale $L^{lc}$ de $L$ et notons $G$ le groupe de Galois $\Gal (L/K)$. Le corps $L^{lc}$ est lui-même une pro-$\ell$-extension de $K$; soit donc $U=\Gal (L^{lc}/K)$ son groupe de Galois et $A$ le sous-groupe $\Gal (L^{lc}/L)$, qui s'identfie au pro-$\ell$-groupe des classes logarithmiques $\Cl_L$. Les sous-groupes de degré nul respectifs de $U$ et de $A$ (\ie les noyaux des morphismes de restriction à $K^c$) sont respectivement $\wi U=\Gal(L^{lc}/K^c)$ et $\wi A= \Gal (L^{lc}/L^c) \simeq \wi\Cl_L$; et $\wi A$ est aussi le sous-groupe dérivé de $U$, puisque $L^c=LK^c$ est, par construction la sous-extension maximale de $L^{lc}$ qui est abélienne sur $K$. L'ensemble de cette discussion peut donc se résumer par le schéma de corps:

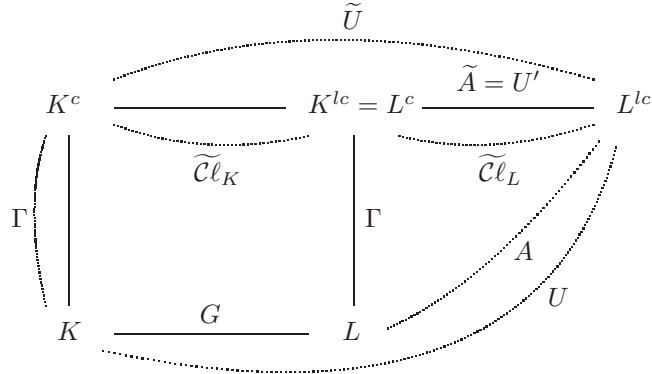
\begin{center}
\unitlength=1.5cm
\begin{picture}(6.6,3)

\put(0.7,0){$K$}
\put(0.8,0.3){\line(0,1){1.5}}
\put(0.6,2){$K^c$}
\put(3.2,0){$L$}
\put(1.2,0.05){\line(1,0){1.7}}
\put(1.95,0.15){$G$}

\bezier{60}(0.6,0.3)(0.4,1.2)(0.6,1.8)
\put(0.3,1.0){$\Gamma$}
\put(3.4,1.0){$\Gamma$}

\put(1.2,2.05){\line(1,0){1.5}}
\put(2.9,2){$K^{lc}=L^c$}
\put(3.9,2.05){\line(1,0){1.5}}
\put(5.6,2){$L^{lc}$}
\put(3.3,0.3){\line(0,1){1.5}}

\bezier{120}(1.2,2.3)(3,3.0)(5.4,2.3)
\put(3.2,2.75){$\wi U$}
%\bezier{50}(3.9,2.1)(4.5,2.3)(5.4,2.1)
\put(4.2,2.2){$\wi A=U'$}

\bezier{50}(1.2,1.9)(2.0,1.6)(2.9,1.8)
\put(1.9,1.4){$\wi\Cl_K$}

\bezier{50}(3.7,1.8)(4.5,1.6)(5.4,1.9)
\put(4.4,1.4){$\wi\Cl_L$}
\bezier{70}(3.6,0.1)(4.5,0.5)(5.45,1.75)
\put(4.7,0.7){$A$}
\bezier{150}(1.1,-0.1)(4.9,-0.9)(5.6,1.7)
\put(5.0,0.3){$U$}

\end{picture}
\end{center}\medskip\bigskip

Maintenant, dans le formalisme du corps de classes, l'extension des idèles (côté arithmétique) correspond au transfert (côté théorie des groupes); de sorte que pour montrer que $\wi\Cl_K$ capitule dans $\wi\Cl_L$, il s'agit de vérifier que le transfert\smallskip

\centerline{$\Ver_{U/A} :U/U' \rightarrow A$}\smallskip

\noindent envoie le sous-groupe $\wi U/U'$ de degré nul de $U/U'$ sur le sous-groupe nul de $\wi A$.\medskip

Il est classique de partir pour cela de l'isomorphisme $G= U/A \simeq \wi U/\wi A$ en faisant choix de représentants des classes de $G$ formé d'éléments de $\wi U$, disons $(u_\tau)_{\tau\in G}$ (avec la convention $u_{\tau^{-1}}=u_\tau^{-1}$) et d'introduire le système de facteurs\smallskip

\centerline{$a_{\sigma,\tau}=u_\sigma u_\tau u_{\sigma\tau}^{-1}\in \wi A$}\smallskip

\noindent dont la classe dans $H^2(G,A)$ définit la loi sur $U$ à isomorphisme près.
\pagebreak

%%%%%%%%%%%%%%%%%%%%%%%%%%%%%%%%%%%%%%%%%%%%%%%%%%%%%%%
\noindent{\bf 4. Interprétation en termes d'algèbre linéaire}\medskip
%%%%%%%%%%%%%%%%%%%%%%%%%%%%%%%%%%%%%%%%%%%%%%%%%%%%%%%

En analogie avec le cas classique, nous pouvons alors constuire un $\Zl$-module résolvant $B$ en formant la somme directe:\smallskip

\centerline{$B= A \oplus I_G = A \oplus (\underset{\tau\ne1}{\oplus}\Zl(\tau-1))$,}\smallskip

\noindent munie de l'action de $G$ définie par:\smallskip

\centerline{$\sigma * a=a^\sigma \qquad \& \qquad \sigma *(\tau-1)= a_{\sigma,\tau} + \sigma(\tau-1)$.}\medskip

\noindent Suivant Artin-Tate (\cf \cite{A-T}, Ch. 13, §4), nous avons un diagramme commutatif:
\begin{displaymath}
\xymatrix{
\Cl_K \ar@{->}[r] ^\sim \ar@{->}[d]_{j_{L/K}} & U/U' \ar@{->}[r] ^\sim \ar@{->}[d]^{\Ver_{U/A}} & B/I_G *B \ar@{->}[d]^{\Tr_{B/A}}\\
\Cl_L \ar@{->}[r] ^\sim  & A \ar@{=}[r] & A
}
\end{displaymath}
où l'isomorphime en haut à droite est donné par le {\em logarithme}: \smallskip

\centerline{$au_\tau U'\mapsto a+(\tau-1)+I_G*B$;}\smallskip

\noindent de sorte que l'homomorphisme d'extension $j_{L/K}$ à gauche correspond au morphisme de transfert $\Ver_{U/A}$ au centre et à la trace $\Tr_{B/A}=\sum_{\tau \in G}\tau $ à droite.\medskip

Comme indiqué dans l'introduction, nous sommes intéressés par la restriction de $j_{L/K}$ au sous-groupe des classes logarithmiques de degré nul. Définissons donc le degré d'un élément $b=a+\sum\lambda_\tau(\tau-1)$ de $B$ par la formule: $\deg b =\deg a$. Le sous-module des éléments de degré nul dans $B$ est alors:\medskip

\centerline{$\wi B= \wi A \oplus I_G =\wi A \oplus (\underset{\tau\ne1}{\oplus}\Zl(\tau-1))$.}\smallskip

\noindent Observons que $\wi B$ contient $I_G*B$; et que le quotient $\wi B/I_G *B$ est encore l'image par le logarithme du sous-groupe $\wi U/U'$ des éléments de degré nul du quotient $U/U'$; ce que nous pouvons résumer par le diagramme commutatif:

\begin{displaymath}
\xymatrix{
\wi\Cl_K \ar@{->}[r] ^\sim \ar@{->}[d]_{j_{L/K}} & \wi U/U' \ar@{->}[r] ^\sim \ar@{->}[d]^{\Ver_{U/A}} & \wi B/I_G *B \ar@{->}[d]^{\Tr_{B/A}}\\
\wi\Cl_L \ar@{->}[r] ^\sim  & \wi A \ar@{=}[r] & \wi A
}
\end{displaymath}

En fin de compte, montrer que les classes logarithmiques de degré nul du corps $K$ capitulent dans $L$ revient à vérifier que pour $B=A \oplus I_G$, muni de l'action de $G$ définie plus haut, l'opérateur trace  $\Tr_{B/A}=\sum_{\tau \in G}\tau $ est nul sur le sous-module $\wi B= \wi A \oplus I_G$. 
\medskip

Avant d'aller plus loin, il peut être intéressant de préciser le dénominateur $I_G*B$ qui intervient dans le quotient à droite du diagramme: l'isomomorphisme canonique de $\Zl[G]$-modules $B/A \simeq I_G$ nous donne l'inclusion: $I_G*B\subset A\oplus I_G^2$.

Pour établir l'égalité, il suffit par conséquent de comparer les deux indices $(\wi B:I_G*B)$ et $(\wi B:\wi A\oplus I_G^2)$. Or, nous avons:
\begin{itemize}
\item d'un côté $\wi B/I_G*B \simeq \wi U/U' \simeq \wi\Cl_K \simeq G$, conformément au diagramme;
\item et directement: $\wi B/(\wi A \oplus I_G^2) = (\wi A \oplus I_G)/(\wi A \oplus I_G^2) \simeq I_G/I_G^2 \simeq G$;
\end{itemize}
d'où l'égalité annoncée:

\centerline{$I_G*B=\wi A \oplus I_G^2$.}
\pagebreak

%%%%%%%%%%%%%%%%%%%%%%%%%%%%%%%%%%%%%%%%%%%%%%%%%%%%%%%
\noindent{\bf 5. Approche par la méthode d'Artin-Furtwängler}\medskip
%%%%%%%%%%%%%%%%%%%%%%%%%%%%%%%%%%%%%%%%%%%%%%%%%%%%%%%

La méthode classique, telle qu'exposée dans \cite{A-T}, consiste à écrire le $\ell$-groupe $G$ comme produit direct de $s$ sous-groupes cycliques $\langle\tau_i\rangle$ d'ordres respectifs $e_i$; et à considérer les éléments $b_i=\tau_i-1$ de $\wi B$  pour $i=1,\cdots,s$.\smallskip

Dans le cas classique le module résolvant $B$ est engendré par les $(b_i)_i$; dans le cadre logarithmique qui nous intéresse ici, c'est plus compliqué, du fait qu'il nous faut considérer le sous-module $\wi B$ des éléments de degré nul.\smallskip 

Relevons $\Gamma=\Gal(K^c/K)$ dans $U=\Gal(L^{lc}/K$ en faisant choix d'un prolongement $\gamma$ à $L^{lc}$ d'un progénérateur arbitraire de $\Gamma$. L'isomorphisme
$$
B/I_G*B \simeq U/U' \simeq \Gamma \times G
$$
montre que $B$ est engendré, comme $\Zl[G]$-module, conjointement par $\gamma$ et les $(b_i)_i$. Et les identités $e_j b_i \in I_G*B$ s'écrivent ainsi:
$$
e_ib_i=\sum_{j=1}^s \lambda_{ij}*b_j +\mu_i *\gamma
$$
avec les $\lambda_{ij}$ et les $\mu_i$ dans $I_G$. Notant alors $M$ la matrice de terme général $m_{ij}=e_i\delta_{ij}-\lambda_{ij}$, nous obtenons:
$$
M *\left[\begin{array}{c}b_1\\ \vdots  \\ b_s \end{array}\right] = \left[\begin{array}{c}\mu_1\\ \vdots  \\ \mu_s \end{array}\right] *\gamma \;;
$$
puis, par multiplication à gauche par la transcomatrice $\wi M$ de $M$:
$$
\det M *\left[\begin{array}{c}b_1\\ \vdots  \\ b_s \end{array}\right] = 
\wi M M *\left[\begin{array}{c}b_1\\ \vdots  \\ b_s \end{array}\right] =
\wi M \left[\begin{array}{c}\mu_1\\ \vdots  \\ \mu_s \end{array}\right] *\gamma = 
\left[\begin{array}{c}\nu_1\\ \vdots  \\ \nu_s \end{array}\right] *\gamma
$$
c'est-à-dire:  $\det M *b_i = \nu_i*\gamma$ pour un $\nu_i$ de $I_G$, pour chaque $i=1,\cdots,s$.\smallskip

L'isomorphisme $\wi B/\wi A\simeq I_G$ montre alors que $\det M \in \Zl[G]$ est un multiple de la trace, disons: $\det M=\kappa \Tr_{B/A}$, pour un $\kappa$ dans $\Zl$.
Et un calcul de degré dans l'algèbre de groupe donne alors directement: \smallskip

\centerline{$\deg (\det M) =  \det ([\deg(m_{ij})]_{ij})=\det ([e_i \delta_{ij}]_{ij})= \prod e_i =\;|\;G\;|= \deg\Tr_{B/A}$;}\smallskip

\noindent d'où, comme attendu: $\det M =\Tr_{B/A}$ et finalement l'inclusion:
$$\Tr_{B/A}(\wi B)\subset I_G*\gamma$$
puisque $\wi B$ est $\Zl[G]$-engendré par les $b_i$ (pour $i=1,\dots,s)$ et $I_G*\gamma$; de sorte que $\Tr_{B/A}(\wi B)$ est engendré par les $\Tr_{B/A}(b_i)$ pour $i=1,\dots,s$. Ainsi:

\begin{Prop}
Avec les notations ci-dessus, l'image dans $\;\Cl_L$ du sous-groupe des classes logarithmiques de degré nul $\;\wi\Cl_K$ s'identifie au sous-module $\Tr_{B/A}(\wi B)$ du $\Zl[G]$-module $ I_G*\gamma$ construit sur un relèvement du groupe procyclique $\Gamma$.
\end{Prop}

\Remarque La méthode ne donne donc pas la trivialité du groupe $\Tr_{B/A}(\wi B)$, comme dans le cas classique, mais seulement une inclusion.

\pagebreak

%%%%%%%%%%%%%%%%%%%%%%%%%%%%%%%%%%%%%%%%%%%%%%%%%%%%%%%
\noindent{\bf 6. Interprétation en termes de théorie des genres}\medskip
%%%%%%%%%%%%%%%%%%%%%%%%%%%%%%%%%%%%%%%%%%%%%%%%%%%%%%%

Pour aller plus loin dans la description de la capitulation logarithmique, nous pouvons nous inspirer des méthodes initiées par T. Tannaka et F. Terada pour généraliser le {\em Théorème de l'idéal principal} dans le contexte de la théorie cyclique des genres (\cf \cite{T-T, Tan, Ter1, Ter2}).\smallskip

Introduisons le sous-corps $L^{lc}_{K^c}$ de $L^{lc}$ fixé par le sous-groupe dérivé $\wi U'$ de $\wi U=\Gal(L^{lc}/K^c)$. Par construction $L^{lc}_{K^c}/K^c$ est la plus grande sous-extension de $L^{lc}/K^c$ qui est abélienne sur $K^c$ et nous avons: $\Gal(L^{lc}_{K^c}/L^c)\simeq\wi A/\wi U'$. Cela étant, $L^c$, qui est la plus grande sous-extension de $L^{lc}_{K^c}/K^c$ qui provient par composition avec $K^c$ d'une extension abélienne de $K$, est donc le corps des genres de $L^{lc}_{K^c}$ relativement à l'extension procyclique $K^c/K$;  en particulier $ \Gal(L^c/K^c)$ est le plus grand quotient de $\wi U/\wi U'=\Gal(L^{lc}_{K^c}/K^c)$ sur lequel $\Gamma=\Gal(K^c/K)$ opère trivialement; et nous avons donc: $\Gal(L^{lc}_{K^c}/L^c)\simeq(\wi U/\wi U')^\omega$; d'où, finalement:
$$
\wi A/\wi U'=(\wi U/\wi U')^{(\gamma-1)},
$$
si $\gamma$ désigne un relèvement arbitraire dans $A=\Gal(L^{lc}/K)$ d'un générateur topologique de $\Gamma$. Et l'ensemble de cette discussion peut être résumé par le schéma de corps:
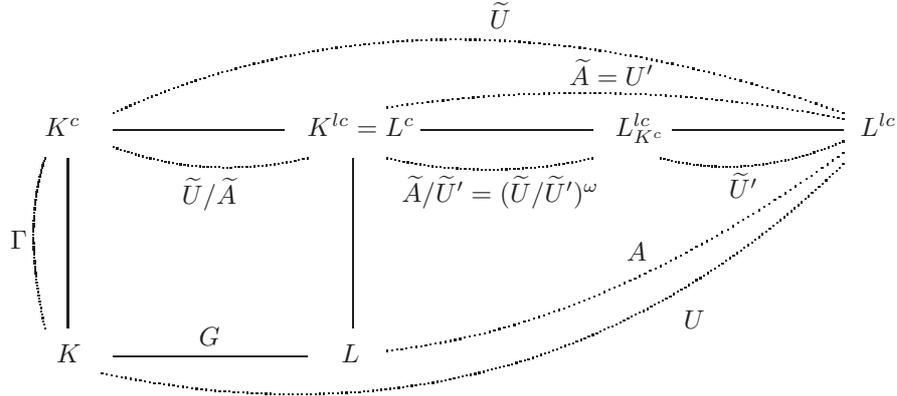
\begin{center}
\unitlength=1.5cm
\begin{picture}(8,3.2)

\put(0.7,0){$K$}
\put(0.8,0.3){\line(0,1){1.5}}
\put(0.6,2){$K^c$}
\put(3.2,0){$L$}
\put(1.2,0.05){\line(1,0){1.7}}
\put(1.95,0.15){$G$}

\bezier{60}(0.6,0.3)(0.4,1.2)(0.6,1.8)
\put(0.3,1.0){$\Gamma$}

\put(1.2,2.05){\line(1,0){1.5}}
\put(2.9,2){$K^{lc}=L^c$}
\put(3.9,2.05){\line(1,0){1.5}}
\put(5.6,2){$L^{lc}_{K^c}$}
\put(3.3,0.3){\line(0,1){1.5}}
\put(7.75,2){$L^{lc}$}
\put(6.1,2.05){\line(1,0){1.5}}

\bezier{150}(1.2,2.2)(4,3.5)(7.6,2.2)
\put(4.5,2.95){$\wi U$}
\bezier{90}(3.6,2.25)(5.6,2.55)(7.6,2.15)
\put(5.2,2.45){$\wi A=U'$}

\bezier{50}(1.2,1.9)(2.0,1.6)(2.9,1.8)
\put(1.8,1.4){$\wi U/\wi A$}

\bezier{50}(3.6,1.8)(4.5,1.6)(5.4,1.8)
\put(3.73,1.42){$\wi A/\wi U'=(\wi U/\wi U')^\omega$}
\bezier{50}(6.0,1.8)(6.8,1.6)(7.6,1.95)
\put(6.6,1.45){$\wi U'$}

\bezier{90}(3.6,0.1)(5.5,0.3)(7.6,1.85)
\put(5.7,0.9){$A$}
\bezier{180}(1.1,-0.1)(4.9,-0.9)(7.6,1.75)
\put(6.2,0.3){$U$}

\end{picture}
\end{center}\medskip\bigskip

\noindent où nous avons noté $\;\omega=\gamma-1$.\medskip

Précisons quelques points.  Comme $\Gamma\simeq\gamma^{\Zl}$ opère trivialement sur le quotient $U/A\simeq G$, nous pouvons écrire:\medskip

\centerline{$u_\tau^\gamma = \gamma u_\tau \gamma^{-1} = a_\tau u_\tau$}\smallskip

\noindent pour tout $\tau\in G$, avec\smallskip

\centerline{$a_\tau = u_\tau^\gamma  u_\tau^{-1} =  [\gamma, u_\tau] =  \gamma^{1-\tau}$;}\medskip

Convenons de noter $\;\wi U^\omega$ le sous-groupe de $\wi A$ qui est engendré par les $a_\tau$ lorsque $\tau$ décrit $G$.  Cela étant, nous avons:

\begin{Lem}
Avec ces conventions, il vient: $\wi A = \wi U'\wi U^\omega$.
\end{Lem}

\Preuve C'est une conséquence immédiate de l'égalité $\wi A/\wi U'=(\wi U/\wi U')^{(\omega)}$ donnée par la théorie des genres.

\pagebreak

%%%%%%%%%%%%%%%%%%%%%%%%%%%%%%%%%%%%%%%%%%%%%%%%%%%%%%%
\noindent{\bf 7. Approche par la méthode de Tannaka-Terada}\medskip
%%%%%%%%%%%%%%%%%%%%%%%%%%%%%%%%%%%%%%%%%%%%%%%%%%%%%%%

Revenons maintenant sur le module résolvant $B$. Rappelons que nous avons fait choix d'un relèvement $\gamma$ dans $A=\Gal(L^{lc}/L)$ d'un générateur topologique de $\Gamma$, ce qui nous a permis de définir les éléments $a_\tau$ pour $\tau \in G$. En termes d'algèbre linéaire, \ie en notations additives, il vient ainsi:\medskip

\centerline{$a_\tau = u_\tau^\gamma  u_\tau^{-1} =  [\gamma, u_\tau] =  \gamma^{1-\tau} = (1-\tau)*\gamma \in I_G*\gamma \subset \wi A$;}\medskip

\noindent et nous pouvons transporter l'action de $\Gamma$ sur $U$ au module $B$ en posant:\medskip

\centerline{$\gamma * (\tau-1) = a_\tau + (\tau -1)$.}\medskip

En particulier $B$ peut ainsi être regardé comme un module sur l'algèbre $\Lambda[G]$ construite sur l'algèbre d'Iwasawa $\Lambda=\Zl[[\omega]]$ en l'indéterminée  $\omega=\gamma-1$. De fait l'identité ci-dessus donne immédiatement:\smallskip

\centerline{$\omega^2*(\tau-1) = \omega*a_\tau =0$,}\smallskip

\noindent de sorte que $B$ est annulé par $\omega^2$.  Plus précisément:

\begin{Lem}
Avec les conventions précédentes, on a: $\omega *\wi  B = \wi  U^\omega =I_G *\gamma$. Et $\wi B$ est $\Lambda[G]$-engendré par les éléments $b_i=\tau_i-1$, pour $i=1,\dots,s$.
\end{Lem}
\medskip

\Preuve Il vient, effet : \smallskip

\centerline{$\omega *\wi  B = \omega *(\wi A\oplus I_G) = \omega *I_G = \underset{\tau \in G}{\sum} \Zl \;\omega *(\tau-1) 
={\sum}_{\tau \in G} \Zl  \;a_\tau =\wi U^\omega$;}\par
\noindent d'où, par passage au quotient à partir de l'isomorphisme: $\wi B/I_G*\wi B \simeq \wi U/\wi U'$:\smallskip

\centerline{$\wi B/(I_G*\wi B +\omega*\wi B) \simeq \wi U/\wi U'\wi U^\omega =  \wi U/\wi A \simeq G$;}\smallskip

\noindent de sorte que $\wi B$, regardé comme $\Lambda[G]$-module, est engendré par les $s$ éléments $b_i=\tau_i-1$ construits sur un système minimal de générateurs de $G$.\medskip

Ce point acquis, les identités $e_j b_i \in I_G*B+\omega *\wi B$ s'écrivent :
$$
e_ib_i=\sum_{j=1}^s \mu_{ij}*b_j +\omega*\sum_{j=1}^s\nu_{ij} *b_j
$$
avec les $\mu_{ij}$ dans $I_G$ et les $\nu_{ij}$ dans $\Zl[G]$. Notant alors $M$ la matrice de terme général $m_{ij}=e_i\delta_{ij}-\lambda_{ij}$ et $N$ celle de terme général $\nu_{ij}$, nous obtenons:
$$
[M-\omega N] *\left[\begin{array}{c}b_1\\ \vdots  \\ b_s \end{array}\right] = \left[\begin{array}{c}0\\ \vdots  \\ 0\end{array}\right] \;;\;
{\rm donc:}\quad \det [M-\omega N] *\wi B = 0.
$$
Un calcul modulo $\wi A$ montre alors, comme précédemment, que le déterminant de la matrice $M$ coïncide avec la trace $\Tr_{B/A}$, de sorte que le déterminant $d= \det [M-\omega N] $ (calculé modulo $\omega^2$) est de la forme:
$$
d=\Tr_{B/A}-\omega\delta
$$
pour un certain $\delta\in\Zl[G]$. En particulier:

\begin{Prop}
Comme opérateur sur $\wi B$, la trace s'écrit: $\Tr_{B/A}=\omega\delta$; et il suit:

\centerline{$\Tr_{B/A}(\wi B)= \delta\omega*\wi B=\delta I_G*\gamma$.}
\end{Prop}

\pagebreak

%%%%%%%%%%%%%%%%%%%%%%%%%%%%%%%%%%%%%%%%%%%%%%%%%%%%%%%
\noindent{\bf 8. Capitulation des classes invariantes}\medskip
%%%%%%%%%%%%%%%%%%%%%%%%%%%%%%%%%%%%%%%%%%%%%%%%%%%%%%%
\medskip

La Proposition 4 ci-dessus précise naturellement l résultat obtenu à la fin de la section 5 par le procédé d'Artin-Furtwängler,  puisqu'elle décrit précisément le sous-module image $\Tr_{B/A}(\wi B)$ des classes logarithmique étendues; mais elle ne donne pas d'information directe sur le sous-groupe de $\wi\Cl_K$ qui capitule dans $\wi\Cl_L$. Nous allons voir que la méthode de Tannaka-Terada donne cependant des éléments de réponse:

\begin{Th}
Le noyau de la trace $\Tr_{B/A}$ regardé dans le quotient $\wi B/I_G*\wi B$ contient le sous-groupe ambige $(\wi B/I_G*\wi B)^\Gamma = {}^{-1}\omega(I_G*\wi B)/I_G*\wi B$.\par
En d'autres termes, on a l'implication:\smallskip

\centerline{$\omega(\wi b) \in I_G*\wi B \Rightarrow  \Tr_{B/A}(\wi b)=0$.}
\end{Th}

\Remarque Le quotient $\wi B/I_G*\wi B\simeq\wi U/\wi U'$ est annulé par $\omega^2$. Le théorème affirme seulement que le sous-groupe annulé par $\omega$ est contenu dans le noyau de $\Tr_{B/A}$.\medskip

\Preuve Le noyau $\wi A \cap I_G*\wi B$ du morphisme naturel $A\rightarrow \wi B/I_G*\wi B \simeq \wi U/\wi U'$ étant $\wi U'$,  nous avons $\wi A \cap I_G*\wi B=\wi U'$, donc:\smallskip

\centerline{$\omega\wi B \cap I_G*\wi B= \omega\wi B \cap \wi A \cap I_G*\wi B = \wi U^\omega \cap \wi U'$.}\smallskip

\noindent Et le Théorème sera établi si nous montrons que l'opérateur $\delta$ est nul sur $\wi U'$.\smallskip

Observons pour cela que le groupe $\wi U$ est engendré par $\wi A$ et les $(u_\tau)_{\tau \in G}$; de sorte que son sous-groupe dérivé $\wi U'$, lui, est engendré conjointement
\begin{itemize}
\item par les commutateurs $[a, u_\tau]=a^{(1-\tau)}$ avec $a\in \wi A$ et $\tau\in G$ et 
\item par les $[u_\sigma, u_\tau]$, lorsque $\sigma$ et $\tau$ parcourent un système générateur de $G$.
\end{itemize}
Introduisons le sous-module $\partial\wi B$ de $\wi B$ qui est engendré par les commutateurs:\smallskip

\centerline{$[b_i,b_j]=b_i*b_j-b_j*b_i=(\tau_i-1)\!*\!(\tau_j-1)\!-\!(\tau_j-1)\!*\!(\tau_i-1)=a_{\tau_i,\tau_j}-a_{\tau_j,\tau_i}$.}
\smallskip

Traduites en notations additives, les observations précédentes nous donnent l'identité:

\centerline{$\wi U' = I_G*A + \partial\wi B$.}\smallskip

\noindent Supposons établie la trivialité de $\delta$ sur le sous-module $\partial\wi B$. Il vient alors :\smallskip

\centerline{$\delta*\wi U' = \delta*I_G*\wi A =\delta I_G*(\wi U'+\omega*\wi B)= I_G*(\delta*\wi U')$;}

\noindent d'où: 

\centerline{$\delta*\wi U'=0$,}\smallskip

\noindent ce qui est précisément le résultat attendu. Ainsi, le Théorème résulte du:

\begin{Lem}[Lemme d'Adachi] 
L'opérateur $\delta$ est nul sur le sous-module $\partial\wi B$.
\end{Lem}

\Preuve Il s'agit de vérifier que l'on a: $\delta b_i*b_j = \delta b_j*b_i$, pour tout $(i,j)$; ce qui se fait en transposant {\em mutatis mutandis} dans le cadre procyclique considéré ici les calculs de déterminants effectués par N. Adachi dans son rapport sur le Théorème de l'idéal principal \cite{Ada} consacré au Théorème de Tannaka-Terada.

\Remarque Le quotient $\wi B/I_G*\wi B\simeq \wi U/\wi U'$ étant réputé fini, le noyau de $\omega$ regardé comme endomorphisme de $\wi B/I_G*\wi B$, a même ordre que son conoyau. On a donc:\smallskip

\centerline{$({}^{-1}\omega(I_G*\wi B):I_G*\wi B)=(\wi U:\wi U'\wi U^\omega)=(\wi U:\wi A) =|G|$.}\smallskip

\noindent Et le {\em noyau de capitulation} est au moins d'ordre $|G|=|\;\wi\Cl_K|$ (mais dans $\wi U/\wi U'$).\medskip
%\pagebreak

%%%%%%%%%%%%%%%%%%%%%%%%%%%%%%%%%%%%%%%%%%%%%%%%%%%%%%%
\noindent{\bf 9. Conclusion}\medskip
%%%%%%%%%%%%%%%%%%%%%%%%%%%%%%%%%%%%%%%%%%%%%%%%%%%%%%%
\medskip

Contrairement à ce qui se produit pour les $\ell$-groupes de classes au sens habituel, on voit que les méthodes classiques, même très élaborées, ne permettent pas de conclure dans le cadre logarithmique à un analogue complet du résultat d'Artin-Furtwängler, puique la capitulation y apparait comme un sous-groupe éventuellment strict du groupe des classes total.\par

De fait cette observation, a priori décevante, se trouve confirmée par l'expérimentation numérique. Donnons un contre-exemple simple qui ruine clairement tout espoir de transposer naïvement dans le cadre logarithmique  le classique théorème 94 de Hilbert sur la capitulation:

\begin{Sco}
Pour $\ell=3$, il existe des  corps quadratiques imaginaires $K$ dont le 3-groupe des classes logarithmiques de degré nul est cyclique d'ordre 3 et qui possèdent une extension cyclique  $L$ de degré $[L:K]=3$  logarithmiquement non ramifiée pour lesquelles l'homomorphisme $j_{L/K}:\wi \Cl_K\rightarrow\wi\Cl_L$ est injectif.
\end{Sco}

\Preuve Sous la conjecture de Gross-Kuz'min (et donc en pratique dans tous les exemples numériques étudiés), le quotient de Herbrand $q(G,\wi\E_L)$ attaché aux unités logarithmiques dans une extension cyclique $L/K$ vaut 1, en vertu de l'expression du caractère du $\Zl[G]$-module $\wi\E_L$ donnée dans \cite{J28}. En particulier, dès que  le groupe $H^2(G,\wi\E_L)$ est trivial, il en est dès lors de même du groupe $H^1(G,\wi\E_L)$ qui s'identifie précisément à la capitulation logarithmique  $\Ker\;j_{L/K}$.\smallskip
Les calculs numériques effectués par Karim Belabas avec {\sc pari}, montrent que c'est, par exemple, le cas pour $K=\Q[\sqrt{-31}]$ et $L=K[x]$, où x est racine du polynôme $X^3 + 3X + (\theta+9)/2$ avec $\theta^2=-31$.\medskip

On voit par là une différence essentielle d'avec le cas des groupes de classes au sens habituel : pour les unités au sens ordinaire, le quotient de Herbrand $q(G,E_L)$ relatif à une extension cyclique de corps de nombres $L/K$ n'est jamais trivial; il est égal au degré $[L:K]$ de l'extension, de sorte que le premier groupe de cohomologie $H^1(G,E_L)$, qui mesure précisément la capitulation dès lors que l'extension considérée $L/K$ est non ramifiée, n'est pas trivial non plus.
Et c'est la clef du théorème 94. Dans le cadre logarithmique en revanche, la trivialité du quotient de Herbrand interdit toute minoration de la capitulation; ce que confirme le calcul.\par

A contrario, le même argument fournit une majoration de la capitulation logarithmique meilleure que celle que l'on a dans le cas classique:

\begin{Prop}
Soient $L/K$ une $\ell$-extension cyclique de corps de nombres logarithmiquement non ramifiée et $G$ son groupe de Galois. Sous la conjecture de Gross-Kuz'min, l'ordre du sous-groupe des classes logarithmiques de $K$ qui capitulent dans $L$ satisfait la majoration:
$$
|\wi{Cap}_{L/K}| = |H^1(G,\wi\E_L)| = |H^2(G,\wi\E_L)| \le (r_K+c_K+\delta_K)^{[L:K]},
$$
où $r_K$ et $c_K$ représentent les nombres de places réelles et complexes du corps $K$ et $\delta_K$ vaut 0 ou 1 suivant que $K$ contient ou pas les racines $\ell$-ièmes de l'unité.
\end{Prop}

\Preuve Le groupe $\wi\E_K$ des unités logarithmiques de $K$ est alors, en effet, le produit du sous-groupe $\mu_K$ des racines d'ordre $\ell$-primaire de l'unité contenues dans $K$ et d'un $\Zl$-module libre de dimension $r_K+c_K$ (cf. \cite{J28}, \S3).

\def\refname{\small{\sc Références}}
\bigskip
{\small
\noindent{\sc Adresse:}
Univ. Bordeaux \& CNRS,
Institut de Mathématiques de Bordeaux, UMR 5251,
351 Cours de la Libération,
F-33405 Talence cedex

\noindent{\sc Courriel:}
 \tt jean-francois.jaulent@math.u-bordeaux1.fr
}

\end{document}